\documentclass[11pt,reqno]{amsart}
\usepackage{amsfonts,amscd,amsmath,amsthm,amssymb,latexsym}
\allowdisplaybreaks[3]
\usepackage{enumerate}

\advance\textwidth by 1.2in \advance\oddsidemargin by -.6in
\advance\evensidemargin by -.6in
\newtheorem*{cor}{Corollary}
\newtheorem*{lem}{Lemma}
\newtheorem*{prop}{Proposition}

\theoremstyle{definition}
\newtheorem*{defn}{Definition}
\theoremstyle{definition}
\newtheorem{thm}{Theorem}

\newtheorem*{rem}{Remark}

\newenvironment{pf}{\proof}{\endproof}
\newcounter{cnt}
\newenvironment{enumerit}{\begin{list}{{\hfill\rm(\roman{cnt})\hfill}}{%
\settowidth{\labelwidth}{{\rm(iv)}}\leftmargin=\labelwidth%
\advance\leftmargin by
\labelsep\rightmargin=0pt\usecounter{cnt}}}{\end{list}}

\theoremstyle{remark}

\numberwithin{equation}{section} 

\def\wh#1{\widehat{#1}}
\def\wt{\operatorname{wt}}

\def\opl_#1^#2{\text{\scriptsize$\bigoplus\limits_{\text{\footnotesize$#1$}}^{\text{\footnotesize$#2$}}$}}

\begin{document}

\newcommand{\thmref}[1]{Theorem~\ref{#1}}
\newcommand{\secref}[1]{Section~\ref{#1}}
\newcommand{\lemref}[1]{Lemma~\ref{#1}}
\newcommand{\propref}[1]{Proposition~\ref{#1}}
\newcommand{\corref}[1]{Corollary~\ref{#1}}
\newcommand{\remref}[1]{Remark~\ref{#1}}
\newcommand{\defref}[1]{Definition~\ref{#1}}
\newcommand{\er}[1]{(\ref{#1})}
\newcommand{\id}{\operatorname{id}}
\newcommand{\tensor}{\otimes}
\newcommand{\nc}{\newcommand}
\newcommand{\rnc}{\renewcommand}
\newcommand{\qbinom}[2]{\genfrac[]{0pt}0{#1}{#2}}
\nc{\cal}{\mathcal} \nc{\goth}{\mathfrak} \rnc{\bold}{\mathbf}
\renewcommand{\frak}{\mathfrak}
\newcommand{\Ext}{\operatorname{Ext}}
\renewcommand{\Bbb}{\mathbb}
\nc\bomega{{\boldsymbol\omega}} \nc\balpha{{\boldsymbol\alpha}}
\nc\bpi{{\boldsymbol\pi}} \nc\bpsi{{\boldsymbol\Psi}}
\nc\bvpi{{\boldsymbol\varpi}}
\newcommand{\lie}[1]{\mathfrak{#1}}
\makeatletter
\def\section{\def\@secnumfont{\mdseries}\@startsection{section}{1}%
  \z@{.7\linespacing\@plus\linespacing}{.5\linespacing}%
  {\normalfont\scshape\centering}}
\def\subsection{\def\@secnumfont{\bfseries}\@startsection{subsection}{2}%
  {\parindent}{.5\linespacing\@plus.7\linespacing}{-.5em}%
  {\normalfont\bfseries}}
\makeatother
\def\subl#1{\subsection{}\label{#1}}
\nc{\Hom}{\operatorname{Hom}}
\nc{\End}{\operatorname{End}}
\nc{\ch}{\operatorname{ch}}
\nc{\ev}{\operatorname{ev}}

\nc{\krsm}{KR^\sigma(m\omega_i)}
\nc{\krsmzero}{KR^\sigma(m_0\omega_i)}
\nc{\krsmone}{KR^\sigma(m_1\omega_i)}
 \nc{\vsim}{v^\sigma_{i,m}}
 \nc{\whg}{\widehat{\lie g}}
 \nc{\Cal}{\cal} \nc{\Xp}[1]{X^+(#1)} \nc{\Xm}[1]{X^-(#1)}
\nc{\on}{\operatorname} \nc{\Z}{{\bold Z}}
\nc{\J}{{\cal J}} \nc{\C}{{\bold C}} \nc{\Q}{{\bold Q}}
\renewcommand{\P}{{\cal P}}
\nc{\N}{{\Bbb N}} \nc\boa{\bold a} \nc\bob{\bold b} \nc\boc{\bold
c} \nc\bod{\bold d} \nc\boe{\bold e} \nc\bof{\bold f}
\nc\bog{\bold g} \nc\boh{\bold h} \nc\boi{\bold i} \nc\boj{\bold
j} \nc\bok{\bold k} \nc\bol{\bold l} \nc\bom{\bold m}
\nc\bon{\bold n} \nc\boo{\bold o} \nc\bop{\bold p} \nc\boq{\bold
q} \nc\bor{\bold r} \nc\bos{\bold s} \nc\bou{\bold u}
\nc\bov{\bold v} \nc\bow{\bold w} \nc\boz{\bold z} \nc\boy{\bold
y} \nc\ba{\bold A} \nc\bb{\bold B} \nc\bc{\bold C} \nc\bd{\bold D}
\nc\be{\bold E} \nc\bg{\bold G} \nc\bh{\bold H} \nc\bi{\bold I}
\nc\bj{\bold J} \nc\bk{\bold K} \nc\bl{\bold L} \nc\bm{\bold M}
\nc\bn{\bold N} \nc\bo{\bold O} \nc\bp{\bold P} \nc\bq{\bold Q}
\nc\br{\bold R} \nc\bs{\bold S} \nc\bt{\bold T} \nc\bu{\bold U}
\nc\bv{\bold V} \nc\bw{\bold W} \nc\bz{\bold Z} \nc\bx{\bold x}
\nc{\whh}{\widehat{\lie h}}
\title[Graded level zero integrable representations]{Graded level zero integrable representations of affine Lie algebras}
\author{Vyjayanthi Chari}
\thanks{This work was partially supported by the NSF grant DMS-0500751}
\address{Department of Mathematics, University of
California, Riverside, CA 92521.}
\email{vyjayanthi.chari@ucr.edu}
\author{Jacob Greenstein}
\email{jacob.greenstein@ucr.edu}
\date{February 22, 2006}
\maketitle
\section{Introduction}
In this paper, we continue the  study (\cite{C,CG,CP}) of the
category $\cal I^{fin}$ of integrable level zero
representations with finite dimensional weight spaces of affine
Lie algebras. The category of integrable representations with
finite dimensional weight spaces and of non-zero level is
semi-simple (\cite{Rao}) and the simple objects are the highest
weight modules. In contrast, it is easy to see that the category
$\cal I^{fin}$ is not semi-simple. For instance, the derived
Lie algebra of the affine algebra itself provides an example of a
non-simple indecomposable object in that category. The simple
objects in~$\cal I^{fin}$ were classified  in~\cite{C,CP} but
not much else is known about the structure of the category.

The category $\cal I^{fin}$ can be regarded as the graded
version of the category $\cal F$ of finite-dimensional
representations of the corresponding loop algebra. The structure
of $\cal F$ has been studied extensively in both the quantum and
classical cases in recent years (\cite{BN,CM,CM1,EM,Nak}). There is a natural functor $L$
from  $\cal F$ to $\cal I^{fin}$ which we study in this paper
(Section~\ref{FuncL}) and which has many nice properties. But
there are  important properties that fail, for example in general
the functor~$L$ only maps irreducible objects to completely
reducible objects. In particular, the trivial representation in
$\cal F$ is mapped to an infinite direct sum of one-dimensional
representations.

The latter phenomenon is  the main source of difficulty and makes
the study of the category $\cal I^{fin}$ interesting in its own
right, since the results and proofs sometimes require substantial
modifications from the  ones for $\cal F$. It is immediate
  from the
preceding comments that objects in the category $\cal I^{fin}$
are not always of finite length. However, in Section~\ref{Fin} of
this paper, we show that a weaker version of the finite length
property holds, namely that the number of non-trivial
constituents of an indecomposable module is finite. Moreover, each
indecomposable module admits an analogue of a composition series
which we call a pseudo-Jordan-H\"older series. Such a series is
unique up to a natural equivalence. As a result we are
able to conclude (\thmref{main1}) that any object in $\cal
I^{fin}$ is a direct sum of indecomposable modules. Moreover,
only finitely many of the  indecomposable summands are non-trivial
modules for the corresponding loop algebra. This, in particular,
implies that $\cal I^{fin}$ has a block decomposition. The
methods used in this part of the paper rely only on facts which
remain valid for quantum affine algebras.

In Section~\ref{Blocks} we obtain a parametrization for the blocks
and describe the blocks explicitly (\thmref{main2}).  The blocks
in~$\mathcal I^{fin}$ are parametrized by the orbits for the natural
action of the group~$\bc^\times$ on the set $\Xi$ of finitely
supported functions from $\bc^\times$ to the quotient group of the
weight lattice of the underlying simple finite dimensional Lie
algebra $\lie g$ by its root lattice. It was proved in \cite{CM}
that $\Xi$ parametrizes blocks in the category~$\mathcal F$. To
prove the result for $\cal I^{fin}$ we use the functor $L$ defined
and studied in Section~\ref{FuncL}. One of the tools used crucially
in \cite{CM} was the notion of the finite-dimensional Weyl module
introduced and studied in \cite{CPweyl} (in the quantum case they
appear in~\cite{K1,K}). One can define in a similar way the notion
of the graded Weyl module, which is an object in $\cal I^{fin}$ and
has the usual universal properties. A natural question  is then
whether the functor $L$ maps Weyl modules in $\cal F$ to graded Weyl
modules or at least to direct sums of graded Weyl modules. This
question turns out to be rather difficult since it is equivalent to
proving the conjecture of~\cite{CPweyl} on the dimension of the Weyl
modules in $\cal F$. This conjecture was established in
\cite{CPweyl} for the affine algebra whose underlying simple Lie
algebra $\lie g$ is isomorphic to~$\lie{sl}_2$. In other cases, as
Hiraku Nakajima has pointed out to us recently, the dimension
conjecture can be deduced as follows. The results of~\cite{K1, K}
imply that the Weyl modules are specializations (the $q=1$ limit)
of certain finite-dimensional quotients of the extremal modules for
the quantum affine algebra. Then it follows from the results in~\cite{BN}
and~\cite{Nak,Nak1} that these quotients (and hence their
specializations) have the correct dimension. Other
approaches to proving this conjecture, which  do not rely on
the quantum case,
have been studied recently: in~\cite{CL} for  the case of $\lie
g$ isomorphic to $\lie{sl}_r$, $r\ge 3$ and in \cite{FoL} for any
simply-laced~$\lie g$.

\subsection*{Acknowledgments}
We are grateful to Hiraku Nakajima for discussions relating to the
dimension conjecture.
\section{Preliminaries}\label{Pre}

Throughout the  paper, $\bz$ (respectively, $\bz_+$, $\bn$) will denote
the set of integers (respectively, non-negative,  positive integers).

\subsection{}\label{P10}
Let $\frak g$ be  a complex finite-dimensional simple Lie algebra
and $\lie h$ a   Cartan subalgebra of~$\lie g$. Set $I=\{1,\dots,
\dim\lie h\}$ and let $\{\alpha_i\,:\,i\in I\}$ (respectively,
$\{\varpi_i:i\in I\}$) be a set of simple roots (respectively, fundamental
weights) of $\frak g$ with respect to $\frak h$, $R^+$ (respectively, $Q$,
$P$) be the corresponding set of positive roots (respectively the root
lattice, the weight lattice). Let $Q^+$, $P^+$ be the
$\bz_+$-span of the simple roots and fundamental weights
respectively.  It is convenient to set $\varpi_0=0$. Let $\ge $ be
the standard partial order on $P$ defined by: $\lambda\ge \mu$ if
$\lambda-\mu\in Q^+$. Let $\theta\in R^+$ be the highest root and
if~$\lie g$ is not simply laced denote by  $\theta_s$ the
highest short root.  Denote by~$W$ the Weyl group of~$\lie g$.

Given $\alpha\in R$, let $\lie g_\alpha$ be the corresponding root
space. For $\alpha\in R^+$, fix non-zero  elements
$x^\pm_\alpha\in\lie g_{\pm\alpha}$ and $\alpha^\vee\in\lie h$,
such that
$$
[\alpha^\vee ,x^\pm_\alpha]=\pm 2x^\pm_\alpha,\qquad
[x^+_\alpha,x^-_\alpha]=\alpha^\vee
$$
and write
$$
\lie n^\pm=\opl_{\alpha\in R^+}^{}\lie
g_{\pm\alpha}=\opl_{\alpha\in R^+}^{}\bc x_{\pm\alpha}.
$$

It is well-known that $\Gamma=P/Q$ is a finite abelian group. For
any $0\ne \gamma\in \Gamma$ fix the unique minimal representative
$0\ne \varpi_\gamma\in P^+$, i.e   if $\lambda\in P^+$ satisfies
$\lambda\le \varpi_\gamma$ then $\lambda=\varpi_\gamma$.
\begin{lem} Let $\lambda\in P^+\setminus\{0\}$ and let $\gamma\in\Gamma$ be such
that $\lambda=\varpi_\gamma\pmod Q$. Then
$\lambda\ge \varpi_\gamma$. In addition, if $\lambda\in
P^+\cap Q^+$ and $\lambda\ne 0$,  then $\lambda\ge\beta$
where $\beta=\theta$ if $\lie g$ is simply laced and
$\beta=\theta_s$ otherwise.
\end{lem}
\begin{pf}
For $\lie g\cong\lie{sl}_{\ell+1}$ one can show by a direct computation
that $\varpi_i+\varpi_j\in \varpi_{i+j\pmod{\ell+1}}+Q^+$, $1\le
i,j\le \ell$. The statement follows since the minimal
representatives of the elements of $\Gamma$ in this case are
$\varpi_i$, $0\le i\le \ell$ and any $\lambda\in P^+\setminus\{0\}$ can be
written as $\lambda=\sum_{r=1}^N \varpi_{i_r}$ for some $N\ge 1$
and $1\le i_r\le \ell$. For all other types the statement is
trivially checked.
\end{pf}

\subsection{}\label{P15}
For $\lambda\in P^+$, let $V(\lambda)$ be the irreducible
finite-dimensional $\lie g$-module with highest weight vector
$v_\lambda$, i.e the cyclic module generated by $v_\lambda$ with
defining relations:
$$
\lie n^+v_\lambda=0,\quad hv_\lambda=\lambda(h)v_\lambda,\quad
(x^-_{\alpha_i})^{\lambda(\alpha_i^\vee)+1}v_\lambda=0,
$$
for all $i\in I$, $h\in\lie h$. Given any $\lie g$-module $M$ and
$\mu\in \lie h^*$ set
$$
M_\mu=\{v\in M: hv=\mu(h)v,\, h\in \lie
h\}.
$$
If $M$ is finite-dimensional then $M=\bigoplus_{\mu\in\lie
h^*}M_\mu,$ and moreover
$$
M\cong\bigoplus_{\lambda\in P^+}V(\lambda)^{\oplus m_\lambda(M)},
\qquad m_{\lambda}(M)\in\bz_+.
$$

The following Lemma is standard (see \cite{H} for instance).
\begin{lem} Let $\lambda, \mu\in P^+$ be such that $\lambda\ge\mu$. Then,
 $V(\lambda)_{\mu}\not=0$.\qed
\end{lem}

\subsection{}\label{P20}
Given  a Lie algebra $\lie a$, let $\bu(\lie a)$ denote  the
universal enveloping algebra of $\lie a$ and let $L(\lie a)=\lie
a\otimes \bc[t, t^{-1}]$  be  the loop algebra  of $\lie a$ with the Lie bracket given
by
$$
[x\tensor f,y\tensor g]=[x,y]_{\lie a}\tensor fg,
$$
for all~$x,y\in\lie a$, $f,g\in\bc[t^{\pm1}]$.  The
Lie algebra $L(\lie a)$ and its universal enveloping algebra  are
$\bz$-graded
by the powers of $t$. We
shall identify $\lie a$ with the subalgebra $\lie a\otimes 1$  of
$L(\lie a)$. Denote by~$L^e(\lie a)=L(\lie a)\oplus\bc d$ the extended loop algebra of~$\lie a$, in which
$[d,x\tensor t^n]=n x\tensor t^n$.

\subsection{}\label{P25}
Let $\lie h_e=\lie h\oplus\bc d$, which is an abelian Lie subalgebra of $L^e(\lie g)$.
Define $\delta\in\lie h_e^*$ by
$$
\delta(d)=1,\qquad \delta(h)=0,\qquad  \forall\, h\in\lie h.
$$

We regard elements of $\lie h^*$ as elements of $\lie h_e^*$ by setting
$\lambda(d)=0$ for $\lambda\in\lie h^*$. In particular, we identify $P$ and
$P^+$ with their respective images in~$\lie h_e^*$. Obviously, $\lie h_e^*=\lie h^*\oplus\bc\delta$.
Set~$P_e=P\oplus\bz\delta$ and let~$P_e^+=P^+\oplus\bz\delta$. Then~$P_e\subset\{\lambda\in\lie h_e^*\,:\,
\lambda(\alpha_i^\vee)\in\bz\,\,\forall\,i\in I\}$.

Let $\wh W$ be the affine Weyl group
associated with~$W$. Its image in the group of automorphisms of~$\lie h_e^*$ identifies
with the semidirect product of $W$ with an abelian group generated by the $t_{h}$, $h \in\bz(W\theta^\vee)\rangle$
where $t_{h} (\lambda)=\lambda-\lambda(h)\delta$ for all $\lambda\in\lie h_e^*$.
Given~$\lambda\in P^+$, let $r_\lambda=\min_{h\in \bz(W\theta^\vee)} \{ \lambda(h)\,:\,
\lambda(h)>0\}\in\bn$.
This number exists since~$\bz(W\theta^\vee)$ is contained in the $\bz$-span of the~$\alpha_i^\vee$, $i\in I$.
Then for all~$s\in\bz$, there exists~$w\in\wh W$ such that $w(\lambda+s\delta)=\lambda+\bar s\delta$,
where $0\le \bar s<r_\lambda$.

\subsection{}\label{P30}  Given $i\in I$, let
$\Lambda^\pm_i(u)\in\bu(L(\lie h))[[u]]$ be defined by
$$\Lambda_i^\pm(u)=\sum_{k\ge 0}\Lambda_{i,\pm
k}u^k=\exp\Big(-\sum_{k>0} \frac{\alpha_i^\vee\otimes t^{\pm
k}}{k}\,u^k\Big).$$
It is easy to see (cf.~\cite{CPweyl}) that $\bu(L(\lie h))$ is the
polynomial algebra on $\Lambda_{i,\pm k}$, $i\in I$, $k\in\bz$,
$k\ne 0$.

\section{Elementary properties of the category $\mathcal I^{fin}$}\label{Elem}

\subsection{}\label{E20}
Recall that a $L^e(\lie g)$-module $V$  is said to be integrable if
$$
V=\bigoplus_{\mu\in\lie h_e^*} V_\mu.
$$
and the elements $x^\pm_\alpha\otimes t^s$ act locally nilpotently
on $V$ for all $\alpha\in R^+$ and $s\in\bz$. Denote by $\mathcal I$ the
category of integrable $L^e(\lie g)$-modules. Let
${\wt}^e(V)=\{\mu\in\lie h_e^*\,:\,V_\mu\ne 0\}$ be the set of
weights of $V$ with respect to $\lie h_e$. It is well-known that the
set ${\wt}^e(V)$ is $\wh W$-invariant.
The following Lemma follows immediately from~\ref{P25} and  will be used repeatedly in the rest of the paper.
\begin{lem} Let $V\in \operatorname{Ob}\cal I$ and assume that
$\mu+s\delta\in {\wt}^e(V)$ where~$\mu\in P^+$, $s\in\bz$. Then $\mu+r\delta\in{\wt}^e(V)$ for
some $0\le r< r_\lambda$.\qed
\end{lem}

\subsection{}\label{E40}
For $a\in \bc$ set
$$
V\{ a\}=\{v\in V\,:\, \text{$dv=(a+k)v$ for some $k\in\bz$} \}.
$$
It is trivial to check that $V\{a\}$ are $L^e(\lie g)$-submodules of $V$ and that
$V\{a\}=V\{b\}$ if and only if $a-b\in\bz$. For any $\overline{a}\in\bc/\bz$
set  $V\{\overline{a}\}=V\{a\}$ where $a$ is any representative of
$\overline{a}$. Let $\mathcal I\{\bar a\}$, $\bar a\in\bc/\bz$ be
the full subcategory of $\mathcal I$ whose objects are
$L^e(\lie g)$-modules $V$ satisfying $V=V\{\bar a\}$.
\begin{lem}
Let $V$ be an integrable $L^e(\lie g)$-module.
\begin{enumerate}[{\rm (i)}]
\item\label{E40.10i} Let $\mu\in{\wt}^e(V)$. Then $\mu(\alpha_i^\vee)\in\bz$,
$i\in I$.

\item\label{E40.10iv}\label{E40.10v} We
have
$$
V=\bigoplus_{\bar a\in\bc/\bz} V\{\bar a\}.
$$
Moreover, if
$V=V\{\bar a\}$, $V'=V'\{\bar b\}$ with $\bar a\not=\bar b$, then
$\Ext^1_{\mathcal I}(V,V')=0$. In particular,
$$
\mathcal I=\bigoplus_{\bar a\in\bc/\bz} \mathcal I\{\bar a\},
$$
and  the categories
$\mathcal I\{\bar a\}$ are equivalent for all $a\in\bc/\bz$.
\end{enumerate}
\end{lem}
\begin{pf}
Part \eqref{E40.10i} is standard and follows from the
representation theory of $\lie{sl}_2$ applied to the subalgebras of~$L^e(\lie g)$
spanned by the elements $\{x_{\alpha_i}^\pm,\alpha_i^\vee\}$ for
$i\in I$. The first two statements
in~\eqref{E40.10iv} are straightforward while for the last it is
sufficient to observe that the functor $V\mapsto V\tensor
\bc_{-a\delta}$, where $\bc_{-a\delta}$ is the 1-dimensional
$L^e(\lie g)$-module on which $L(\lie g)$ acts trivially and $d$ acts
by $-a$, provides an equivalence of categories between $\mathcal
I\{\bar a\}$ and $\mathcal I\{\bar 0\}$.
\end{pf}
It follows, in particular, that we can restrict ourselves to the subcategory~$\mathcal I\{\bar0\}$
of~$\mathcal I$. Observe that $V\in\operatorname{Ob}\mathcal I$ is an object in~$\mathcal I\{\bar0\}$ if and only if~$\wt^e(V)\subset P_e$.
Let $\mathcal I^{fin}$ be the subcategory of $\cal I$ consisting
of modules $V$ such that $\wt^e(V)\subset P_e$ and $\dim V_\mu<\infty$ for all $\mu\in
P_e$.

\subsection{}\label{E60}

Let $V$ be an integrable $L^e(\lie g)$-module. For $\lambda\in P^+_e\setminus\bz\delta$,
set $V^+_\lambda=\{v\in V_\lambda: L(\lie n^+)v=0\}$ and
$$
V^+=\bigoplus_{\lambda\in P_e^+\setminus\bz\delta\,:\, 0\le \lambda(d)< r_\lambda} V_\lambda^+.
$$
\begin{prop}
Let $V$ be an object in $\mathcal I^{fin}$ and suppose that  $V_1\subsetneq
V_2\subsetneq\cdots$ is  an ascending chain of
$L^e(\lie g)$-submodules of $V$. Then ${\wt}^e(V_r/V_{r-1})\subset
\bz\delta$ for all but finitely many $r\ge 1$.  In particular,
$V^+$ is finite-dimensional.
\end{prop}
\begin{pf}
Write
$$
V=\bigoplus_{\gamma\in \Gamma} V[\gamma],
$$
where
$$
V[\gamma]=\displaystyle\bigoplus_{\mu=\varpi_\gamma\!\pmod
Q,\,n\in\bz} V_{\mu+n\delta}.
$$
Since $V[\gamma]$ is obviously a $L^e(\lie g)$-submodule of
$V$ and $\Gamma$ is a finite group we can assume without
loss of generality  that $V=V[\gamma]$. Suppose for a
contradiction that ${\wt}^e(V_r/V_{r-1})\not\subset \bz\delta$
for infinitely many $r\ge 1$. Since~$V$ is integrable, $V_r/V_{r-1}$ is a (possibly infinite)
direct sum of finite dimensional $\lie g$-modules. In particular, for infinitely many
$r\ge 1$ there exists $\mu_r\in P^+\setminus\{0\}$ and $s_r\in\bz$ such that~$v\in (V_r/V_{r-1})_{\mu_r+s_r\delta}$,
$v\not=0$ generates a simple highest weight $\lie g$-submodule isomorphic to~$V(\mu_r)$.
By~\lemref{P10}, $\mu_r\ge\varpi_\gamma$, hence  $V(\mu_r)_{\varpi_\gamma}\not=0$ by~\lemref{P15} and we conclude that
$\varpi_\gamma+s_r\delta\in{\wt}^e(V_r/V_{r-1})$.

Suppose first that
$\varpi_\gamma\ne 0$. Then by~\lemref{E20}, $\varpi_\gamma+
\bar s_r\delta\in{\wt}^e(V_r/V_{r-1})$ with~$0\le \bar s_r < r_{\varpi_\gamma}$.
It follows that there exists $0\le
s< r_{\varpi_\gamma}$  such that
$\varpi_\gamma+s\delta\in\wt^e(V_r/V_{r-1})$ for infinitely many~$r\ge1$. This implies that $V_{\varpi_\gamma+s\delta}$ is infinite dimensional
which is clearly a contradiction. If $\varpi_\gamma=0$ and
$\mu_r\ne 0$, then by~\lemref{P10} we have that
$\mu_r\ge\beta $ where~$\beta=\theta$ or~$\beta=\theta_s$. The  preceding argument would then imply
that $V_{\beta+s\delta}$ is infinite dimensional for some~$0\le s<r_\beta$, which is again a contradiction.
\end{pf}

\subsection{}\label{E80}
In the category~$\mathcal I^{fin}$ the role similar to that of highest weight modules is played by
$\ell$-highest weight modules.
\begin{defn} Let $V$ be an integrable module. A non-zero element $v\in V^+_\lambda$ is called an
$\ell$-highest weight vector if~$\bu(L^e(\lie h))v$ is an indecomposable
$\bu(L^e(\lie h))$-submodule of~$V$ and if
$\dim(\bu(L^e(\lie h))v)_{\lambda+r\delta}{\le 1}$ for all~$r\in\bz$. We say
that $V$ is $\ell$-highest weight if $V=\bu(L^e(\lie g))v$ where $v\in V$
is  an $\ell$-highest weight vector.
\end{defn}
It is not hard to
see by the usual arguments that an  $\ell$-highest weight module
is indecomposable. Observe also that if~$v\in V^+_\lambda$ and~$\bu(L^e(\lie h))v$ is a simple $\bu(L^e(\lie h))$-module, then
$\dim(\bu(L^e(\lie h))v)_{\lambda+r\delta}\le 1$ and so~$v$ is an~$\ell$-highest weight vector.

\subsection{}\label{E100}
Retain the assumption of~\ref{E80} and
consider $V^{L(\lie g)}=\{v\in V\,:\, L(\lie g)v=0\}$. Obviously, $V^{L(\lie g)}$ is
an $L^e(\lie g)$-submodule of $V$.
\begin{lem}\begin{enumerate}[{\rm(i)}]
\item\label{E100.i} $V^{L(\lie g)}\cong\bigoplus_{r\in\bz} \bc_{r\delta}^{m_r}$, $m_r\in\bz_+$
\item\label{E100.ii} Suppose that ${\wt}^e(V)\subset \bz\delta$. Then~$V=V^{L(\lie g)}$.
\item\label{E100.iii} Suppose that $V\not=V^{L(\lie g)}$. Then $V/V^{L(\lie g)}$ does not admit an $L^e(\lie g)$-submodule isomorphic to
$\bc_{r\delta}$, $r\in\bz$.
\end{enumerate}
\end{lem}
\begin{pf}
Parts~\eqref{E100.i} and \eqref{E100.ii} are immediate. For~\eqref{E100.iii},
let $V_0=V^{L(\lie g)}\subsetneq V$ and suppose that there exists $v\in V\setminus V_0$ such that $L(\lie g)v\subset V_0$.
Since $V$ is a weight module, we can write, uniquely, $v=\sum_{k} v_{\mu_k}$, $\mu_k\in\lie h_e^*$. Since $hv\in V_0$ for all
$h\in\lie h_e$ and
${\wt}^e(V_0)\subset\bz\delta$ by~\eqref{E100.i}, it follows that $\mu_k\in\bz\delta$. Then $(x^\pm_\alpha\tensor t^n)v=0$
and so $L(\lie g)v=0$. Thus, $v\in V_0$, which is a contradiction.
\end{pf}

\section{Finiteness in the category $\mathcal I^{fin}$}\label{Fin}

The main result of this section is the following theorem. Although
we state and prove this result only for the classical affine Lie
algebras, it is clear that the proof goes over verbatim to the
quantum case.
\begin{thm}\label{main1}
Let $V$ be an object in $\mathcal I^{fin}$.
Then $V$ is isomorphic
to a direct sum of indecomposable modules. More precisely,   there
exists submodules $U_j$, $j=1,2$ such that $V=U_1\oplus U_2$ with
$U_1\subset V^{L(\lie g)}$, $U_2^+=V^+$.  Moreover, $U_2$ is
isomorphic to a finite direct sum of indecomposable modules.
\end{thm}
\noindent We prove this result in the rest of the section.

\begin{rem}
Note that there exist indecomposable modules of infinite length in $\cal I^{fin}$.
For example, let $M=\lie g\oplus\bc$ as a $\lie g$-module with the
$L(\lie g)$-module structure defined by $(x\tensor t^k)(y\oplus a)=[x,y]\oplus k (x,y)_{\lie g}$, where
$(\cdot\,,\,\cdot)_{\lie g}$ is the Killing form of~$\lie g$,
for all $x,y\in\lie g$, $k\in\bz$ and~$a\in\bc$.
Then $M$ is an indecomposable $L(\lie g)$-module.
Applying the functor~$L$ (cf.~\ref{L60}) we conclude that $L(M)$ contains
$L(\bc)=\bigoplus_{r\in\bz} \bc_{r\delta}$ as a submodule and $L(M)/L(\bc)\cong L(\lie g)$ which is
a simple $L^e(\lie g)$-module. Thus, $L(M)$ has infinite length and it is easy to see that $L(M)$ is indecomposable.
\end{rem}

\subsection{}\label{F30}
The following proposition was established in~\cite{CG} (Proposition~3.2) in the quantum case.
We provide its proof here for the sake of completeness.
\begin{prop}
Let $V$ be an object in~$\mathcal I^{fin}$.
Then there exists
$\lambda\in \wt^e(V)$
such that $\lambda+\eta\notin \wt^e(V)$ for all $\eta\in Q^+\setminus\{0\}$. In particular,
$\lambda\in P^+_e$ and $V^+_\lambda\not=0$.
\end{prop}
\begin{pf}
By~\lemref{E100}, we may assume that~$\wt^e(V)\not\subset\bz\delta$.
Suppose that for each~$\mu\in\wt^e(V)$ there exists~$\nu\in
Q^+$ such that~$\mu+\nu\in\wt^e(V)$. Fix some~$\mu\in\wt^e(V)$. Then
there exists an infinite sequence $\{\eta_r\}_{r\ge 1}$ such
that~$\eta_r \le \eta_{r+1}$ and $\mu+\eta_r\in\wt^e(V)$ for all~$r\ge1$.
Set $W_r:=\bu(\lie g)V_{\mu+\eta_r}$. Then~$W_r$ is an integrable
$\lie g$-module with finite-dimensional weight spaces
and hence is isomorphic to a finite direct sum of simple finite dimensional
$\lie g$-modules $V(\mu_{r,s})$ for some~$\mu_{r,s}\in P^+$.
Choose $s_1$ such that $\nu_1:=\mu_{1,s_1}>\mu$. Such~$s_1$ exists since~$\mu+\eta_1$ is
a weight of~$W_1$. Furthermore, let $r_2$ be the smallest positive
integer so that there exists $s_2$ with
$\nu_2:=\mu_{r_2,s_2}>\mu$, $\mu_{1,s_1}\ne \mu_{r_2,s_2}$. Notice that $r_2$ always exists since the
module $W_1$ is finite-dimensional and the maximal weights which
occur in $W_r$ keep increasing. Repeating this process, we obtain
an infinite collection of elements $\nu_k>\mu$, $k\ge 1$ such that
such that $V(\nu_k)$ is isomorphic to an irreducible $\lie
g$-submodule $W(\nu_k)$ of $V$. By~\lemref{P15} it follows
that $V_{\mu}\cap W(\nu_k)\ne 0$ for all~$k\ge1$. Since all
the~$\nu_k$ are distinct, the sum of~$W(\nu_k)$ is direct, which
contradicts the finite-dimensionality of~$V_\mu$. In
particular~$\lambda+\alpha_i\notin \wt^e(V)$ for all~$i\in I$ which
implies that $\lambda\in P^+_e$.
\end{pf}

\subsection{}\label{F40}
\begin{prop}
Suppose that $V$ is an object in $\mathcal I^{fin}$ and assume that $\wt^e(V)$ is
not a subset of $\bz\delta$.  Then $V$ contains an
$\ell$-highest weight submodule generated by~$v\in V^+_{\lambda+r\delta}$ for some $\lambda\in P^+\setminus\{0\}$ and $0\le r<
r_\lambda$.
\end{prop}
\begin{pf}
By~\propref{F30} and~\lemref{E20} we find that there exists
$\lambda\in P^+\setminus\{0\}$ and $0\le s< r_\lambda$ such that~$V_{\lambda+s\delta}^+\not=0$.
Fix a non-zero~$v\in V^+_{\lambda+s\delta}$. Let $V_1=\bu(L^e(\lie g))v$.
We claim that $V_1$ and hence $V$ contains an $\ell$-highest weight
$L^e(\lie g)$-submodule. For this it suffices to prove that the
$\bu(L^e(\lie h))$-module $\bigoplus_{r\in\bz} (V_1)_{\lambda+r\delta}$
contains an irreducible $\bu(L^e(\lie h))$-submodule. If not, then there
exists integers $r_k\in\bz$, and non-zero elements $v_k\in
(V_1)_{\lambda+r_k\delta}$, $k\in\bn$,  such that $\bu(L^e(\lie h))
v_k\supsetneq \bu(L^e(\lie h))v_{k+1}$ and we can assume as usual that
$0\le r_k<r_\lambda$ for all $k$. Since the $v_k$ are obviously
linearly independent, this contradicts the fact that the weight
spaces $V_{\lambda+s\delta}$ are  finite dimensional. Thus, there
exists an element $v_0\in V^+_\mu\setminus\{0\}$ for some~$\mu\in \wt^e(V)$ such that
$\bu(L^e(\lie h))v_0$ is a simple $\bu(L^e(\lie h))$-module. Then $\bu(L^e(\lie g))v_0$ is an $\ell$-highest weight module.
\end{pf}
\begin{cor} Suppose that $V\in\operatorname{Ob}\cal I^{fin}$.  Then $V$ contains a simple
$\ell$-highest weight submodule.
\end{cor}
\begin{pf} By~\lemref{E100}, we may assume that $\wt^e(V)$ is not a subset of $\bz\delta$.
By the above Proposition, it is enough to consider the case when $V$ is an $\ell$-highest
weight module generated by~$v\in V^+_{\lambda+k\delta}$, $\lambda\in P^+$, $0\le k<r_\lambda$.
Assume for a contradiction that the corollary
is false. Then there exists an infinite family of non-zero
$\ell$-highest weight $L^e(\lie g)$-submodules $V_r$, $r>0$ of $V$, such
that $V_r\supsetneq V_{r+1}$. Let $v_r$ be an $\ell$-highest
weight vector generating $V_r$ and assume that the weight of $v_r$ is
$\mu_r+s_r\delta$, where $\mu_r\in P^+\cap (\lambda-Q^+)$,
$s_r\in\bz$. As before, we may also assume that $0\le s_r< r_{\mu_r}$. Since the set  $ P^+\cap(\lambda-Q^+)$ is
finite, the weights of the linearly independent elements $v_r$,
$r\ge 1$  are contained in a finite set, which is a contradiction since
weight spaces of $V$ are finite dimensional.
\end{pf}

\subsection{}\label{F60}
The next proposition provides an analogue of the finite length property for the category~$\mathcal I^{fin}$.
\begin{defn}
We say that a descending chain of $L^e(\lie g)$-submodules $V=V_N\supsetneq V_{N-1} \supsetneq \cdots \supsetneq V_1\supsetneq V_0=0$ of a
$L^e(\lie g)$-module $V\in\operatorname{Ob}\mathcal I^{fin}$ is a {\em pseudo-Jordan-H\"older series}, if
\begin{enumerate}[(i)]
\item\label{F60.i} either $V_k/V_{k-1}$ is a simple $\ell$-highest weight module and ${\wt}^e(V_k/V_{k-1})\not\subset
\bz\delta$
\item\label{F60.ii} or $V_k/V_{k-1}\cong\bigoplus_{r\in\bz} \bc_{r\delta}^{m_r}$ for some $m_r\in\bz_+$
\end{enumerate}
and for all~$k$, at most one of $V_k/V_{k-1}$, $V_{k+1}/V_k$ satisfy~\eqref{F60.ii}.
\end{defn}
\begin{prop}
Let $V\in \operatorname{Ob}\mathcal I^{fin}$ be  indecomposable. Then~$V$ admits a pseudo-Jordan-H\"older series.
In particular, $V$ has finite length, in the usual sense, if and only if the set of $r\in\bz$
such that $\bc_{r\delta}$ is a subquotient  of $V$ is finite.
\end{prop}
\begin{pf} Since $V$ is indecomposable ${\wt}^e(V)\not\subset\bz\delta$. If
$V^{L(\lie g)}\not=0$, set $V_1=V^{L(\lie g)}$. Otherwise, by
Corollary~\ref{F40}, choose  $V_1$ to be a simple $\ell$-highest weight
submodule  generated by an element of $V_{\lambda_1+r\delta}^+$
for some $\lambda_1\in P^+\setminus\{0\}$ and $0\le r<r_\lambda$.
This procedure can obviously be repeated to get an ascending chain satisfying
the conditions of the above definition.
Proposition~\ref{E60} implies that we must reach a stage when
$(V/V_N)^{L(\lie g)}=V/V_N $ and hence $V_{N+1}=V$ which proves
the finiteness condition of the proposition.
\end{pf}

\subsection{}\label{F80}
Let~$V$ be an object in~$\mathcal I^{fin}$ and let~$V=V_N\supseteq \cdots \supseteq V_1\supseteq V_0=0$ and
$V=V_{N'}'\supseteq \cdots \supseteq V_1'\supseteq V_0'=0$
be descending chains of its $L^e(\lie g)$-submodules. We call these
chains equivalent if
\begin{enumerate}[(i)]
\item $|\{0< i\le N\,:\, {\wt}^e(V_i/V_{i-1})\not\subset\bz\delta\}|=|\{0< j\le N'\,:\, {\wt}^e(V'_j/V_{j-1}')\not\subset\bz\delta\}|$.
\item For each~$0< i\le N$ such that~${\wt}^e(V_i/V_{i-1})\not\subset\bz\delta$, there exists~$0< j\le N'$ such that
$V_i/V_{i-1}\cong V'_j/V_{j-1}'$.
\end{enumerate}
It is easy to see that the above relation is indeed an equivalence.
\begin{prop}
Let $V$ be an indecomposable object in~$\mathcal I^{fin}$. Then its pseudo-Jordan-H\"older series is unique up to equivalence
defined above.
\end{prop}
\begin{pf}
It is easy to see that any refinement of a pseudo-Jordan-H\"older series
for an indecomposable object~$V$ in~$\mathcal I^{fin}$ is equivalent, in the above sense, to that pseudo-Jordan-H\"older series
itself. The statement follows immediately from the Schreier Refinement Theorem.
\end{pf}
It follows from the above proposition that the number
$$
|\{0< i\le N\,:\, {\wt}^e(V_i/V_{i-1})\not\subset\bz\delta\}|,
$$
where~$V=V_N\supsetneq\cdots \supsetneq V_1\supsetneq V_0=0$ is a pseudo-Jordan-H\"older series for~$V$, is well-defined. We call this number the pseudo-length of~$V$.

\subsection{}\label{F100}
We can now prove~\thmref{main1}.
\begin{pf}
Let $V\in \operatorname{Ob}\mathcal I^{fin}$. For each $s\in\bz$, let
$m_s\in\bz_+$ be the multiplicity of $\bc_{s\delta}$ as a direct
summand of $V$ and set $U_1=\bigoplus_{r\in\bz}
\bc^{m_r}_{r\delta}$.  Clearly $U_1\subset V^{L(\lie g)}$. For
$s\in\bz$ with  $m_s\not=0$ fix a $L^e(\lie g)$-module complement $V^{(s)}$
of $\bc_{s\delta}^{m_s}$  in $V$ and let
$U_2=\bigcap_{r\in\bz\,:\, m_r\not=0} V^{(r)}$. Since $V^{(r)}\cap
\bc^{m_r}_{r\delta}=0$ it follows that $U_1\cap U_2=\{0\}$.

Let $v\in V_\mu$,
$\mu\in{\wt}^e(V)$. If $\mu\not=s\delta$ for some $s\in\bz$ such
that $m_s\not=0$, then $v\in V^{(s)}$ for all $s\in\bz$ and so $v\in
U_2$. Otherwise, since $V=V^{(s)}\oplus \bc_{s\delta}^{m_s}$, we can
write $v$  uniquely as  $v=u_1+u_2$ for some
$u_1\in\bc_{s\delta}^{m_s}$ and $u_2\in V^{(s)}$.
It follows from weight
considerations  that
$u_2\in U_2$ and so we have proved that  $V=U_1\oplus
U_2$. Note that if $U_2\ne 0$, then by~\lemref{E100}
${\wt}^e(U_2)$ is not a subset of $\bz\delta$.

It remains to show that $U_2$ is a finite direct sum of indecomposable
$L^e(\lie g)$-modules. If $U_2$ is indecomposable we are done.
Otherwise, we can write $U_2= M_1\oplus M_2$. Note that ${\wt}^e(M_j)$
is not a subset of $\bz\delta$ for $j=1,2$ since otherwise by
\lemref{E100} we would have a contradiction to the definition
of $U_1$. Hence $\dim M_j^+>0$ for $j=1,2$ and since~$U_2^+=M_1^+\oplus M_2^+$, $\dim U_2^+=\dim
M_1^++\dim M_2^+$. The statement follows by repeating the argument
and using the fact that $\dim U_2^+<\infty$ (cf.~\propref{E60}).
\end{pf}

Since by~\thmref{main1} every object in the category $\mathcal
I^{fin}$ is a direct sum of indecomposables, it makes sense to
define the blocks in that category.

\subsection{}\label{F120}
Let $\mathcal C$ be an abelian category in which any object is a direct sum of indecomposables.
We say that two indecomposable objects $U_i$, $i=1,2$ in an
$\mathcal C$ are linked and write $U_1\sim U_2$ if there do not exist
full abelian subcategories $\mathcal C_i$, $i=1,2$ such that
$U_i\in\operatorname{Ob}\mathcal C_i$ and $\mathcal C= \mathcal
C_1\oplus\mathcal C_2$. If the $U_i$ are decomposable, then they
are said to be linked if every indecomposable summand of $U_1$ is
linked to every indecomposable summand of $U_2$. This defines an
equivalence relation on $\mathcal C$ and  the equivalence classes
are called blocks. Each block is a full abelian subcategory and
the category $\cal C$ is a direct sum of the blocks.

It is not hard to see that the following is an equivalent definition of linking.
Two indecomposable objects $U$, $V$ in~$\mathcal C$ are
linked if and only if there exists a family of indecomposable objects
$U_1=U, U_2,\dots, U_l=V$ in $\mathcal C$ such that either
$\Hom_{\mathcal C}(U_k,U_{k+1})\not=0$ or $\Hom_{\mathcal
C}(U_{k+1},U_k)\not=0$ for all $k=1,\dots,l-1$.

\subsection{}\label{F140}
Following \cite{CM}, let $\Xi$ be the set  of functions
$\chi:\bc^\times\to \Gamma$ with finite support. Addition of
functions defines on ~$\Xi$ the structure of an abelian group.
Given $i\in I$, $a\in\bc^\times$, set $\chi_{i,a}(z)=\delta_{a,z}
\bar\varpi_i$, where $\bar\varpi_i$ denotes the canonical image of
$\varpi_i$ in $\Gamma$. Clearly $\Xi$ is the free abelian group
generated by the $\chi_{i,a}$, $i\in I$, $a\in\bc^\times$.

Define an action of  $\bc^\times$ on $\Xi$ by
$$(a\cdot\chi)(z):=\chi(az),\qquad a, z\in\bc^\times.$$ Let
$\overline{\Xi}$ be the set of orbits in $\Xi$ for this action and
let  $\bar\chi$ the $\bc^\times$-orbit of $\chi\in\Xi.$  In the
rest of the paper, we shall prove the following:
\begin{thm}\label{main2} Assume that $\lie g$ is simply-laced.
The blocks in the category $\mathcal I^{fin}$ are parametrized
by the elements of $\overline{\Xi}$.
\end{thm}

\section{The category $\cal F$ and the functor $L$}\label{FuncL}

\subsection{}\label{L10} Let $\mathcal P^+$ be the set of $I$-tuples of
polynomials $\bpi=(\pi_i)_{i\in I}$ in $u$ with constant term $1$.
We regard $\mathcal P^+$ as a commutative monoid with
multiplication defined component-wise. Let~$\boldsymbol
1=(1,\dots,1)$ and, for $i\in I$ and $a\in\bc^\times$, let
$\bvpi_{i,a}\in\cal P^+$ be the $I$-tuple of polynomials with the
polynomial  $(1-au)$ in the $i^{th}$ place and one everywhere
else. The elements $\bvpi_{i,a}$, $i\in I$ and $a\in\bc^\times$
generate the monoid $\cal P^+$. For $\bpi\in \cal P^+$, set
$\lambda_\bpi=\sum_{i\in I} (\deg\pi_i)\varpi_i\in P^+$.
Conversely, given $\lambda\in P^+$ and $a\in\bc^\times$, set
$\bpi_{\lambda,a}=\prod_{i\in I}\bvpi_{i,a}^{\lambda(\alpha_i^\vee)}$.

We say that $\bpi=(\pi_i)_{i\in I},\bpi'=(\pi_i')_{i\in
I}\in\mathcal P^+$ are co-prime if for all~$i,j\in I$, the
polynomials $\pi_i$ and $\pi'_j$ are co-prime.
Clearly any
$\bpi\in\cal P^+$ can be written, uniquely, as a product
$\bpi_{\lambda_j,a_j}$ for some $\lambda_j\in P^+$ and
$a_j\in\bc^\times$, $1\le j\le k$ and  $a_r\ne a_j$ if $r\ne j$.

Let $m:\cal P^+ \to \bz_+$ be the map defined by
setting  $m(\bpi)$ to be  the maximal non-negative integer~$r$
such that $\pi_i\in \bc[u^r]$ for all $i\in I$.
\begin{lem}
Suppose  that $\lambda,\mu\in P^+$ and $\lambda>\mu$.
Let $\bpi\in \mathcal P^+$ and $a\in\bc^\times $ be such that
$\bvpi_{i,a}$ is coprime to $\bpi$ for some $i\in I$. Then either
$m(\bpi\bpi_{\lambda,a})=1$ or $m(\bpi\bpi_{\mu,a})=1$ .
\end{lem}
\begin{pf}
Observe that for all~$\bpi'=(\pi_i')_{i\in I}\in\mathcal P^+$, $m(\bpi')>1$ implies that~$\sum_{r} \mu_r(h) b_r=0$
for all~$h\in\lie h$, provided that~$\bpi'=\prod_r \bpi_{\mu_r,b_r}$ where the~$b_r$ are distinct and~$\mu_r\in P^+$.
Indeed, for~$i\in I$, $\sum_r \mu_r(\alpha_i^\vee)
b_r$ is obviously the coefficient of~$u$ in the polynomial~$\pi'_i(u)$. Since we
assume that~$m(\bpi')>1$, it follows from the
definition of~$m(\bpi')$ that the coefficient of~$u$
in all the~$\pi'_i$, $i\in I$ must be zero.

Write~$\bpi=\prod_{r=1}^k \bpi_{\lambda_r,a_r}$ where the~$a_r$ are distinct.
Obviously, $a\not=a_r$ for all~$1\le r\le k$.
Let~$\beta=\lambda-\mu\in Q^+\setminus\{0\}$. Suppose that
$m(\bpi_\lambda\bpi) > 1$. By the above argument, $\lambda(h)a+\sum_r \lambda_r(h) a_r=0$ for all~$h\in \lie h$.
Since~$\beta\not=0$, there exists~$i\in I$ such that~$\beta(\alpha_i^\vee)\not=0$. Then
$\mu(\alpha_i^\vee)a+\sum_r \lambda_r(\alpha_i^\vee) a_r=
(\lambda-\beta)(\alpha_i^\vee) a+\sum_r \lambda_r(\alpha_i^\vee)a_r=-\beta(\alpha_i^\vee)a\not=0$
which implies~${m(\bpi_\mu\bpi)=1}$.
\end{pf}

\subsection{}\label{L20}
\begin{defn} Let $V$ be an $L(\lie g)$-module and let
$0\ne v\in V$. We say that $v$ is an $\ell$-highest vector
if~$$L(\lie n^+)v=0,\ \ \dim_\bc\bu(L(\lie h))v=1.$$ The module~$V$
is said to be $\ell$-highest weight if it is generated, as $L(\lie
g)$-module, by an $\ell$-highest weight vector.\end{defn}
Given
$\bpi=(\pi_i)_{i\in I}\in P^+$, let $\mathcal W(\bpi)$ be the $L(\lie g)$-module generated by
a vector $w_\bpi$ satisfying the relations:
$$
L(\lie n^+)w_\bpi=0,\qquad hw_\bpi=\lambda_\bpi(h)w_\bpi,\qquad
\Lambda_{i,\pm k}w_\bpi=\pi^\pm_{i,k}w_\bpi,
$$ where $h\in\lie h$,
$\pi_i(u)=\sum_{k\ge 0}\pi_{i,k}^+ u^k$ and $\sum_{k\ge 0} \pi_{i,k}^- u^k=
u^{\deg\pi_i}\pi_i(u^{-1})/(u^{\deg\pi_i}\pi_i(u^{-1})|_{u=0})$. The modules
$\mathcal W(\bpi)$ are clearly $\ell$-highest weight modules for $L(\lie
g)$.

Let $\cal F$ be the category of finite-dimensional representations
of $L(\lie g)$.
The following was proved in~\cite{C,CP,CPweyl}.
\begin{prop}\label{wpi}
\begin{enumerate}[\rm(i)]
\item\label{wpi.i} For all $\bpi\in\cal P^+$,  $\mathcal W(\bpi)$ is an indecomposable object in $\mathcal F$.
\item\label{wpi.ii} Any $\ell$-highest weight module in $\cal F$ is a quotient of
$\mathcal W(\bpi)$ for some $\bpi\in\cal P^+$.
\item\label{wpi.iii} The modules $\mathcal W(\bpi)$
have a unique irreducible quotient $\mathcal V(\bpi)$ and any irreducible
module in $\cal F$ is isomorphic to $\mathcal V(\bpi)$ for some $\bpi\in\cal
P^+$.
\item\label{wpi.iv} Let $\bpi,\bpi'\in\cal P^+$ be coprime. Then
$\mathcal W(\bpi)\otimes \mathcal W(\bpi)\cong \mathcal W(\bpi\bpi')$. In particular, any
quotient of $\mathcal W(\bpi\bpi')$ is isomorphic to a tensor product of
quotients of $\mathcal W(\bpi)$ and $\mathcal W(\bpi')$.
\item\label{wpi.v} Let~$a\in\bc^\times$ and let~$\tau_a$ be the automorphism of~$L(\lie g)$ defined by~$\tau_a(x\tensor t^r)=
a^r x\tensor t^r$ for all~$x\in\lie g$, $r\in\bz$. Then $\tau_a^* \mathcal W(\bpi)$ is isomorphic, as
an $L(\lie g)$-module, to $\mathcal W(\bpi(au))$. \qed
\end{enumerate}
\end{prop}
\begin{cor}
Let~$\bpi\in\cal P^+$ and suppose that~$m=m(\bpi)>1$. Then there exists~$\eta_m\in\End_{\bc}\mathcal W(\bpi)$ such
that $\eta_m(w_\bpi)=w_\bpi$ and $\eta_m( (x\tensor t^n)w)=\zeta^{-r}_m (x\tensor t^r) \eta_{\bpi}(w)$ for all~$x\in\lie g$, $r\in\bz$ and $w\in \mathcal W(\bpi)$,
where~$\zeta_m$ is an $m$th primitive complex root of unity. Moreover, $\eta_m$ is of order~$m$.
\end{cor}
\begin{pf}
Let~$a\in\bc^\times$. Since~$\tau_a^*\mathcal W(\bpi)$ is isomorphic to~$\mathcal W(\bpi)$ as a vector space, it follows from~\eqref{wpi.v} that there exists a
map~$\eta_{\bpi,a}\in\Hom_{\bc}(\mathcal W(\bpi),\mathcal W(\bpi(au)))$ such that~$\eta_{\bpi,a}(w_\bpi)=w_\bpi$ and
$\eta_{\bpi,a}( (x\tensor t^r) w)=a^r (x\tensor t^r) \eta_{\bpi,a}(w)$ for all~$x\in\lie g$, $r\in\bz$ and $w\in W$.
Set~$\eta_m=\eta_{\bpi,\zeta_m^{-1}}$. Since~$m=m(\bpi)>1$, $\bpi(\zeta_m u)=\bpi(u)$ and so~$\eta_m\in\End_{\bc}\mathcal W(\bpi)$.
The first two properties of~$\eta_m$ are immediate. For the last, observe that since~$\mathcal W(\bpi)$ is generated by~$w_\bpi$
as an $L(\lie g)$-module, it follows from the properties of~$\eta_m$ that $\mathcal W(\bpi)$ is a direct sum of eigenspaces of~$\eta_m$
all eigenvalues of~$\eta_m$ are $m$th complex roots of unity and $\zeta_m$ is an eigenvalue of~$\eta_m$.
\end{pf}

\begin{rem}
The module~$\mathcal V(\bpi_{\lambda,a})$ is isomorphic to~$V(\lambda)$ as a $\lie g$-module, the $L(\lie g)$-module
structure being defined by the evaluation at~$a$, that is, $(x\tensor t^n) v=a^n xv$ for all~$x\in\lie g$,
$n\in\bz$ and~$v\in V(\lambda)$.
More generally, $\mathcal V(\bpi)$, $\bpi\in\cal P^+$ is isomorphic to a tensor product of modules of the form~$\mathcal V(\bpi_{\lambda,a})$ with
distinct~$a$.
\end{rem}

\subsection{}\label{L40}
The assignment $\bvpi_{i,a}\to \chi_{i,a}$ extends to
a surjective map of monoids $\cal P^+\to\Xi$. Denote by~$\chi_\bpi$ the image of~$\bpi\in \mathcal P^+$ under this map.
Given $\chi\in\Xi$, let $\cal F_\chi$ be the full
subcategory of $\cal F$ whose objects have the following property: $\mathcal V(\bpi)$ is
an irreducible constituent of $V$ in~$\mathcal F_\chi$ only if $\chi=\chi_\bpi$. The
following result was proved in \cite{CM}.
\begin{thm}
We have $\cal F=\bigoplus_{\chi\in\Xi}\cal F_\chi$. Moreover the
$\cal F_\chi$ are the blocks in $\cal F$. \qed
\end{thm}

\subsection{}\label{L60} Define a functor $L:\cal F\to \cal I^{fin}$ by
$L(V)=V\otimes \bc[t,t^{-1}]$ with the $L^e(\lie g)$-module
structure given by:
$$
(x\tensor t^k)(v\tensor t^n)=(x\tensor t^k)v\tensor t^{k+n},\quad
d(v\tensor t^n)=nv\tensor t^n
$$
for all~$x\in\lie g$, $v\in V$ and $k,n\in\bz$. Clearly $L$
preserves direct sums and short exact sequences.

It should be noted that the functor~$L$ is not essentially surjective. For example, the indecomposable
module $L(\lie g)\oplus\bc$ with the~$L^e(\lie g)$-module structure given by
$$
(x\tensor t^k)(y\tensor t^r,a)=([x,y]\tensor t^{r+k}, r\delta_{r,-k}(x,y)_{\lie g}),\quad d(y\tensor t^r,a)=(r y\tensor t^r,0)
$$
is not isomorphic to~$L(V)$ for any~$V$ in~$\cal F$.

Given any $\ell$-highest weight $L(\lie g)$-module generated by an $\ell$-highest weight vector
$v$, let $L^s(V)$ be the $L^e(\lie g)$-submodule of $L(V)$ generated by
$v\otimes t^s$. The following result was proved in \cite{CP} (see also \cite{CG}).
\begin{prop} Let $\bpi\in\cal P^+$, $\bpi\not=\boldsymbol 1$.
\begin{enumerate}[{\rm(i)}]
\item For $0\le s<m(\bpi)$ the module
$L^s(\mathcal V(\bpi))$ is an irreducible $L^e(\lie g)$-submodule of
$L(\mathcal V(\bpi))$ and moreover, $$L(\mathcal V(\bpi))=\bigoplus_{s=0}^{m(\bpi)-1}
L^s(\mathcal V(\bpi)).$$
Further, as $L(\lie g)$-modules we have,
$L^s(\mathcal V(\bpi))\cong L^r(\mathcal V(\bpi))$ for all $0\le s,r<m(\bpi)$.
\item
Any irreducible object in $\cal I^{fin}$ is isomorphic to
$L^s(\mathcal V(\bpi))$ for some $\bpi\in\cal P^+$, $\bpi\not=\boldsymbol1$ and $0\le s<m(\bpi)$ or
to $\bc_{r\delta}$ for some~$r\in\bz$.
\item As $L^e(\lie g)$-modules $L^s(\mathcal V(\bpi))\cong L^r(\mathcal V(\bpi'))$
if and only if~$\bpi'(u)=\bpi(au)$ for some~$a\in\bc^\times$ and~$r=s\pmod{m(\bpi)}$.\qed
\end{enumerate}
\end{prop}

\subsection{}\label{L80}
Motivated by the preceding result, we define an action of
the  group $\bc^\times$  on $\mathcal P^+$ by $(a\cdot
\bpi)(u)=\bpi(au)$ for all $a\in\bc^\times$, $\bpi\in\mathcal P^+$
and we define $\overline{ \cal P^+}$ and $\overline \bpi$ in the
obvious way.  It is easily
checked that the surjective map of monoids $\mathcal P^+\to\Xi$ described in~\ref{L40} induces a surjective map $\overline{\cal P^+}\to
\overline{ \Xi}$. Denote by $\overline{\chi_\bpi}$ the image under
this map of the element $\overline \bpi$.

Given $\overline{\chi}\in\overline{\Xi}$ with
$\overline{\chi}\not=0$, let $\cal I^{fin}_{\overline \chi}$ be
the full subcategory of $\cal I^{fin}$ consisting of
$L^e(\lie g)$-modules $V\in\operatorname{Ob} \cal I^{fin}$ satisfying:
$L^s(\mathcal V(\bpi))$ is an irreducible constituent of $V$ only if
$\overline\chi=\overline{\chi_\bpi}$. We say that  $V$ is in $\cal
I^{fin}_{\overline 0}$ if all its simple constituents are either
of the form $\bc_{r\delta}$, $r\in\bz$ or of the form
$L^s(\mathcal V(\bpi))$  with $\chi_\bpi=0$.

\begin{prop}
\begin{enumerate}[{\rm(i)}]
\item\label{L80.i} If $V\in\operatorname{Ob}\cal F_\chi$ then $L(V)\in\operatorname{Ob}\cal
I^{fin}_{\overline\chi}$.
\item\label{L80.ii}
Suppose that $V$ is an $\ell$-highest weight quotient of $\mathcal W(\bpi)$ for
some~$\bpi\in\mathcal P^+$ and let~$v$ be the canonical image of~$w_\bpi$ in~$V$.
For $0\le s<m(\bpi)$ the submodule  $L^s(V)=\bu(L^e(\lie g))(v\tensor t^s)\subset L(V)$ is an $\ell$-highest weight $L^e(\lie g)$-module
and $$L(V)=\bigoplus_{s=0}^{m(\bpi)-1}L^s(V).$$
\end{enumerate}
\end{prop}
\begin{pf}
Since~$V$ is finite dimensional, it has a Jordan-H\"older series
$V=V_N\supsetneq \cdots \supsetneq V_0=0$. By applying the
functor~$L$ we obtain a filtration $L(V)=L(V_N)\supsetneq \cdots
\supsetneq L(V_0)=0$. Observe that~$L(V_i)/L(V_{i-1})\cong
L(V_i/V_{i-1})$ and thus is either simple or completely reducible.
In the second case, it is a finite direct sum of simple objects
unless $V_i/V_{i-1}\cong\bc$ in which case
$L(V_i/V_{i-1})\cong\bigoplus_{r\in\bz} \bc_{r\delta}$. Therefore,
refining and dropping terms if necessary we can construct a
pseudo-Jordan-H\"older series for~$L(V)$ out of Jordan-H\"older
series for~$V$. The statement follows since the
pseudo-Jordan-H\"older series is unique up to equivalence defined
in~\ref{F80}.

It suffices to prove~\eqref{L80.ii} when~$V=\mathcal W(\bpi)$ and when
$m=m(\bpi)>1$.  Clearly,
$$
L(\mathcal W(\bpi))= \sum_{s\in\bz} L^s(\mathcal W(\bpi)),\qquad
L^s(\mathcal W(\bpi))=L^r(\mathcal W(\bpi))\quad \text{if $s=r\pmod m$}.
$$
Define~$\wh\eta_m\in\End_{\bc} L(\mathcal W(\bpi))$
by $\wh\eta_m(w\tensor t^r)=\zeta_m^r \eta_m(w)\tensor t^r$, where~$\zeta_m$ is an $m$th primitive root of unity and
$\eta_m$ is a map from~\corref{wpi} (cf.~\cite[2.6]{CG}).
It is easy to check that~$\wh\eta_m\in\End_{L^e(\lie g)} L(\mathcal W(\bpi))$ and since $\wh\eta_m$ is clearly of order~$m$, it defines a representation of~$\bz/m\bz$ on~$L(\mathcal W(\bpi))$. It follows that~$L(\mathcal W(\bpi))$ is a direct sum of~$\bz/m\bz$-isotypical components corresponding
to $m$ distinct irreducible characters of the finite abelian group~$\bz/m\bz$. It remains to
observe that~$\bz/m\bz$ acts by its irreducible character corresponding to $\zeta_m^s$ on~$L^s(\mathcal W(\bpi))$ and hence~$L^s(\mathcal W(\bpi))$ is contained in a $\bz/m\bz$-isotypical component of~$L(\mathcal W(\bpi))$.
\end{pf}

\subsection{}\label{L90}
For any object $V$ in $\mathcal I^{fin}$, let  $V^\#=\bigoplus_{\mu\in P_e} V^*_\mu\subset V^*$
be the graded dual of~$V$. Then $V^\#$ is an  $L^e(\lie g)$-submodule of
$V^*$ and is in $\cal I^{fin}$ and the functor sending $V$ to $V^\#$ is exact and contravariant.
It is easy to see that $L(\mathcal V(\bpi))^\# \cong L(\mathcal V(\bpi^*))$, with $\bpi^*\in\mathcal P^+$
satisfying $\lambda_{\bpi^*}=-w_\circ \lambda_\bpi$ where~$w_\circ$ is the longest element of~$W$.
In particular, if~$\bpi_j\in\mathcal P^+$, $j=1,2$, and $\lambda_{\bpi_1}-\lambda_{\bpi_2}\in Q^+\setminus\{0\}$, then
$\lambda_{\bpi^*_2}-\lambda_{\bpi^*_1}\notin Q^+$.

\subsection{}\label{L100}
\begin{prop}
Let~$a\in\bc^\times$ and let~$\lambda_j\in P^+$, $j=1,2$ be such that
$\lambda_1>\lambda_2$. Set~$\bpi_j=\bpi_{\lambda_j,a}$ and suppose that $\Ext^1_{\mathcal
F}(\mathcal V(\bpi_1),\mathcal V(\bpi_2))\not=0$.  Assume
that  $\bpi\in\cal P^+$ is coprime to $\bvpi_{i,a}$ for $i\in I$.
\begin{enumerate}[\rm(i)]
\item\label{L100.ii} If $m(\bpi_1\bpi)=1$ then
$L(\mathcal V(\bpi_1\bpi))$ and $L(\mathcal V(\bpi_2\bpi))$ are linked.
\item\label{L100.iii} If $m(\bpi_1\bpi)>1$ then
$L(\mathcal V(\bpi_2\bpi))$ is simple and is linked to $L^s(\mathcal V(\bpi_1\bpi))$ for
some~$0\le s<m(\bpi_1\bpi)$.
\end{enumerate}
\end{prop}
\begin{pf}
Suppose that we have a non-split short exact
sequence of finite dimensional $L(\lie g)$-modules
$$
0\longrightarrow \mathcal V(\bpi_2)\longrightarrow V\longrightarrow \mathcal V(\bpi_1)\longrightarrow 0.
$$
Since tensoring with $\mathcal V(\bpi)$ is exact and the
functor $L$ preserves short exact sequences, we have a short exact
sequence of $L^e(\lie g)$-modules,
\begin{equation}\label{L100.10}
0\longrightarrow  L(\mathcal V(\bpi_2)\tensor \mathcal V(\bpi))\longrightarrow  L(V\otimes \mathcal V(\bpi))\longrightarrow  L(\mathcal V(\bpi_1)\tensor \mathcal V(\bpi))\longrightarrow  0.
\end{equation}
It is not hard to see (cf.~\cite{CM}) that the module $V$ is an
$\ell$-highest weight module for $L(\lie g)$ and is a quotient
of~$\mathcal W(\bpi_1)$. Moreover,  since $\bpi$ and $\bpi_j$, $j=1,2$ are
coprime, it follows from Proposition \ref{wpi} that $V\otimes
\mathcal V(\bpi)$ is an $\ell$-highest weight  quotient of~$\mathcal W(\bpi_1\bpi)$
and also that  $\mathcal V(\bpi_j)\otimes \mathcal V(\bpi)\cong \mathcal V(\bpi_j\bpi)$,
$j=1,2$.

If~$m(\bpi_1\bpi)=1$ then by~\propref{L80}\eqref{L80.ii} we see
that~$L(V\tensor \mathcal V(\bpi))$ is indecomposable and part~\eqref{L100.ii} is
immediate. Otherwise, by~\lemref{L10}, $m(\bpi_2\bpi)=1$ and so~$L(\mathcal V(\bpi_2\bpi))$ is simple.
By~\propref{L80}\eqref{L80.ii},
$L(V\tensor \mathcal V(\bpi))=\bigoplus_{s=0}^{m(\bpi_1\bpi)-1}
L^s(V\tensor \mathcal V(\bpi))$ and
$L(\mathcal V(\bpi_1\bpi))=\bigoplus_{s=0}^{m(\bpi_1\bpi)-1}
L^s(\mathcal V(\bpi_1\bpi))$. It follows
from~\eqref{L100.10} that the sequence
$$
\begin{CD}
0 \longrightarrow  L(\mathcal V(\bpi_2\bpi))\longrightarrow  L^s(V\otimes \mathcal V(\bpi))\longrightarrow  L^s(\mathcal V(\bpi_1\bpi))\longrightarrow  0
\end{CD}
$$
is exact for some~$0\le s<m(\bpi_1\bpi)$. Part~\eqref{L100.iii} follows since~$L^s(V\tensor \mathcal V(\bpi))$
is $\ell$-highest weight and hence indecomposable.
\end{pf}

\begin{cor} Let $\bpi\in\cal P^+$ and assume that $m(\bpi)> 1$.
 Let  $a\in\bc^\times$ be such that $\bvpi_{i,a}$, $i\in I$ is coprime to $\bpi$.
Then~$L(\mathcal V(\bpi\bpi_{\theta,a}))\sim L(\mathcal V(\bpi))$. In particular,
the modules $L^s(\mathcal V(\bpi))$ and $L^r(\mathcal V(\bpi))$ are linked for any
$0\le s, r<m(\bpi)$.
\end{cor}\label{rslink}
\begin{pf}
By~\lemref{L10}, $m(\bpi_{\theta,a}\bpi)=1$. Since $\Ext^1_{\mathcal F}(\mathcal V(\bpi_{\theta,a}),\bc)\ne 0$ (cf.~\cite{CM}),
the result is immediate from the proposition.
\end{pf}

\subsection{}\label{L120}
\begin{prop} Let $\bpi,\bpi'\in\cal P^+$ and suppose that~$\chi_{\bpi'}=0$. Then
$L(\mathcal V(\bpi\bpi'))\sim L(\mathcal V(\bpi))$.
\end{prop}
\begin{pf}
Clearly it is sufficient to prove the statement for~$\bpi'=\bpi_{\beta,a}$ where
$\beta\in P^+\cap Q\subset Q^+$.

Suppose first that~$\bpi_{\beta,a}$ is co-prime with~$\bpi$.
Since~$\lambda\in P^+\cap Q^+$, by~\cite[Proposition~1.2]{CM} there exists a sequence $\gamma_r\in
P^+\cap Q^+$, $r=0,\dots,N$ such that~$\gamma_0=\beta$, $\gamma_N=0$ and
$\Ext^1_{\mathcal F}(\mathcal V(\bpi_{\gamma_r,a}),\mathcal V(\bpi_{\gamma_{r+1},a}))\not=0$. It follows then
from~\propref{L100} and its corollary that the module
$L(\mathcal V(\bpi_{\gamma_r,a}\bpi))$ is  linked  to
$L(\mathcal V(\bpi_{\gamma_s,a}\bpi))$ for all~$0\le r,s\le N$, which implies the assertion.

If~$\bpi_{\lambda,a}$ is not co-prime with~$\bpi$, when we can write~$\bpi=\bpi_{\mu,a}\bpi_1$
where~$\bpi_1$ is co-prime with~$\bvpi_{i,a}$, $i\in I$. Then~$\bpi\bpi_{\beta,a}=\bpi_{\mu+\beta,a}\bpi_1$.
Again, by~\cite[Proposition~1.2]{CM}, there exists a sequence~$\nu_r\in P^+$, $r=0,\dots,K$ such that~$\nu_0=
\mu+\beta$, $\nu_K=\mu$ and $\Ext^1_{\mathcal F}(\mathcal V(\bpi_{\nu_r,a}),\mathcal V(\bpi_{\nu_{r+1},a}))\not=0$.
Then it follows from~\propref{L100} and its corollary that $L(\mathcal V(\bpi_{\nu_r,a}\bpi_1))\sim L(\mathcal V(\bpi_{\nu_s,a}\bpi_1))$
for all~$0\le r,s\le K$.
\end{pf}

\begin{cor}
Let~$\bpi_j\in\mathcal P^+$, $j=1,2$ be such that
$\overline{\chi}_{\bpi_1}=\overline{\chi}_{\bpi_2}$.  Then
$L(\mathcal V(\bpi_1))\sim L(\mathcal V(\bpi_2))$.
\end{cor}
\begin{pf}
Suppose first that~$\overline{\chi}_{\bpi_j}=0$, $j=1,2$. Then~$\chi_{\bpi_j}=0$ and
it follows from the Proposition
that the $L(\mathcal V(\bpi_j))$, $j=1,2$ are linked to~$L(\bc)$ and hence are linked.

Suppose that~$\overline{\chi_{\bpi_j}}\not=0$. We may assume, without loss of generality, that $\chi_{\bpi_1}=
\chi_{\bpi_2}$. By the Proposition, we may assume that~$\bpi_j=\prod_{r=1}^k \bpi_{\lambda_{j,r},a_r}$
such that~$\lambda_{1,r}-\lambda_{2,r}=\beta_r^+-\beta_r^-$, $\beta_r^\pm\in Q^+$.
Applying the Proposition again we conclude that
$L(\mathcal V(\bpi_2))\sim L(\mathcal V(\bpi_2\prod_{r=1}^k \bpi_{\beta_r^+,a_r}))=L(\mathcal V(\bpi_1\prod_{r=1}^k \bpi_{\beta_r^-,a_r}))
\sim L(\mathcal V(\bpi_1))$.
\end{pf}

\section{Characters of indecomposable objects in~$\mathcal I^{fin}$}\label{Blocks}

Throughout this section we assume that~$\lie g$ is simply laced.

\subsection{}\label{CI20}
In order to complete the proof of~\thmref{main2}, it remains to establish the following proposition.
\begin{prop}
Let~$V$ be an indecomposable object in~$\mathcal I^{fin}$. Then~$V$ is an object in~$\mathcal I^{fin}_{\bar\chi}$
for some~$\chi\in\Xi$.
\end{prop}
\noindent
In order to prove this proposition we will need to establish some properties of $\ell$-highest weight $L^e(\lie g)$-modules.

\subsection{}\label{CI40} Given $\lambda\in P^+_e$, let
$W^e(\lambda)$ be the left $L^e(\lie g)$-module generated by an element
$w_\lambda$ subject to the relations
$$
L(\lie n^+)w_\lambda =0,\quad hw_\lambda=\lambda(h)w_\lambda,\quad
(x_{\alpha}^-\otimes 1)^{\lambda(\alpha^\vee)+1}w_\lambda=0,
$$
for all $\alpha\in R^+$ and $h\in\lie h_e $. Clearly, right
multiplication by elements of $\bu(L(\lie h))$ defines a structure
of a right $\bu(L(\lie h))$-module on $W^e(\lambda)$. The following
can be found in \cite{CPweyl}.
\begin{prop}
\begin{enumerit}
\item Let $J_\lambda=\operatorname{Ann}_{\bu(L(\lie h))} w_\lambda$. Given
$\bpi\in\cal P^+$ with $\lambda_\bpi=\lambda$,  there exists a
maximal ideal $J_\bpi\supset J_\lambda$ in $\bu(L(\lie h))$ such
that one has an isomorphism of $L(\lie g)$-modules
$$
\mathcal W(\bpi)\cong W^e(\lambda)\otimes_{\bu(L(\lie h))}\bu(L(\lie h))/J_\bpi.
$$
Conversely if $J$ is any maximal ideal
in $\bu(L(\lie h))/J_\lambda$ then $J=J_\bpi$ for some
$\bpi\in\cal P^+$ and
$$
\mathcal W(\bpi)\cong
W^e(\lambda)\otimes_{\bu(L(\lie h))}\bu(L(\lie h))/J.
$$
\end{enumerit}
\end{prop}

\subsection{}\label{CI50}
As explained in the introduction, the next theorem
follows from~\cite{BN}, \cite{K,K1}. It was also proved, by
different methods, in \cite{CPweyl} for $\lie g=\lie{sl}_2$, in
\cite{CL} for $\lie g=\lie{sl}_n$, $n\ge 3$ and in \cite{FoL} for
$\lie g$ of all simply laced types.
\begin{thm} \label{free}
Let~$\lambda\in P^+$.
Then there exists~$N_\lambda>0$ such that~$\dim \mathcal W(\bpi)=N_\lambda$ for all~$\bpi\in\mathcal P^+$
satisfying~$\lambda_\bpi=\lambda$.
\end{thm}

\begin{cor}
As a $\bu(L(\lie h))/J_\lambda$-module, $W^e(\lambda)$ is free of rank $N_\lambda$.
\end{cor}
\begin{pf}
Let~$S_\lambda=\bu(L(\lie h))/J_{\lambda}$. Then any maximal ideal in~$S_\lambda$ is
of the form~$J_{\bpi}$ for some~$\bpi\in\cal P^+$ with~$\lambda_\bpi=\lambda$.
Since~$\cal W(\bpi)\cong W^e(\lambda)\tensor_{S_\lambda}S_\lambda/J_{\bpi}$, it follows
from the Theorem that~$\dim W^e(\lambda)\tensor_{S_\lambda} S_\lambda/J=N_\lambda$
for any maximal ideal~$J$ in~$S_\lambda$. The statement follows by Nakayama's Lemma.
\end{pf}

\subsection{}\label{CI80}
Let $K_{\bpi}$ be the kernel of the homomorphism $\psi_\bpi:\bu(L(\lie
h))\to\bc[t^{\pm1}]$ of $\bz$-graded algebras defined by
$\psi_{\bpi}(\Lambda_{i,\pm s})=\pi_{i,s}^\pm t^{\pm s}$. Then
$K_\bpi\supset J_{\lambda_\bpi}$ and we set
$$W^e(\bpi)=W^e(\lambda_\bpi)\tensor_{\bu(L(\lie h))} \bu(L(\lie
h))/K_\bpi.$$
It is clear that any $\ell$-highest weight
$L^e(\lie g)$-module is a quotient of $W^e(\bpi)\tensor \bc_{s\delta}$ for some
$\bpi\in\cal P^+$, $s\in\bz$.
\begin{prop} The module $W^e(\bpi)$ is a free module over $\bu(L(\lie
h))/K_\bpi$ of rank $\dim \mathcal W(\bpi)$. In particular $W^e(\bpi)\cong
L^0(\mathcal W(\bpi))$ as $L^e(\lie g)$-modules.
\end{prop}
\begin{pf}
The first statement is a standard consequence of the corollary of
Theorem~\ref{free}.
For the second, let~$m=m(\bpi)$ and let~$\eta_m\in\End_{\bc}\mathcal W(\bpi)$ be the map from~\corref{wpi}.
By the proof of~\propref{L80}\eqref{L80.ii},
$L^0(\mathcal W(\bpi))$ is spanned by elements~$w\tensor t^s$, $s\in\bz$, where~$w\in \mathcal W(\bpi)$ satisfies~$\eta_m(w)=\zeta_m^{-s} w$,
$\zeta_m$ being an $m$th complex primitive root of unity.
Since $\pi_{i,s}^\pm=0$ unless~$m$ divides~$s$, it follows that the formula
$$
(w\otimes t^r)\Lambda_{i,\pm s}=\pi_{i,s}^\pm (w\otimes t^{r\pm s}),
$$
defines a right action of~$\bu(L(\lie h))$ on~$L^0(\mathcal W(\bpi))$ which
commutes with the left $L(\lie g)$-module action.
It follows that $L^0(\mathcal W(\bpi))$ is a
free module for $\bu(L(\lie h))/K_\bpi$ of rank $\dim(\mathcal W(\bpi))$.

The proposition follows by noticing that  $L^0(\mathcal W(\bpi))$ is an
$L^e(\lie g)$-module quotient of $W^e(\bpi)$, that the quotient
map is also a map of right $\bu(L(\lie h))$-modules and that~$L^0(\mathcal W(\bpi))$ and~$W^e(\bpi)$ are
free $\bu(L(\lie h))/K_\bpi$-modules of the same rank.
\end{pf}

\begin{cor}
Let~$V$ be an $\ell$-highest weight $L^e(\lie g)$-module.
Then $V$ is a quotient of $L^s(\mathcal W(\bpi))$ for some $\bpi\in\mathcal{P}^+$, $0\le s<m(\bpi)$.
\end{cor}

\subsection{}\label{CI100}
The crucial point in the proof of~\propref{CI20} is the following result.
\begin{prop}
Suppose that $V_i$, $i=1,2$ are simple objects in~$\mathcal I^{fin}_{\bar\chi_i}$ respectively,
with $\bar\chi_1\not=\bar\chi_2$. Then~$\Ext^1_{\mathcal I^{fin}}(V_1,V_2)=0$.
\end{prop}
\begin{pf}
Suppose
that there exists   a non-split short exact sequence of $L^e(\lie
g)$-modules
$$
\begin{CD}
0 @>>> V_2 @>\phi>> M @>\xi>> V_1 @>>>0.
\end{CD}
$$
Since~$\bc_{r\delta}\in \operatorname{Ob}\mathcal I^{fin}_{\bar 0}$ for all~$r\in\bz$, we may assume
that at least one of $V_1$ or $V_2$ is not isomorphic to $\bc_{r\delta}$ for $r\in\bz$. By
passing to graded duals if necessary (cf.~\ref{L90}), we can assume that $V_1\cong L^s(\mathcal V(\bpi_1))$ for  some $\bpi_1\in\cal P^+$ with $\bpi_1\ne\mathbf 1$,
$0\le s<m(\bpi_1)$ and also that if $V_2\cong L^r(\mathcal V(\bpi_2))$, $0\le r<m(\bpi_2)$ then $\lambda_{\bpi_2}-\lambda_{\bpi_1}\notin Q^+$.
Since $L^s(\mathcal V(\bpi_1))\cong
L^0(\mathcal V(\bpi_1))\tensor \bc_{s\delta}$, we further assume by tensoring the above
short exact sequence, if necessary, with~$\bc_{-s\delta}$
that $s=0$. Choose $m_{\bpi_1}\in M$ such that $\xi(m_{\bpi_1})=v_{\bpi_1}\tensor 1$.
Since the sequence is not split and  $\lambda_{\bpi_2}\not > \lambda_{\bpi_1}$, we have
$$
L(\lie n^+)m_{\bpi_1}=0,\qquad M=\bu(L^e(\lie g))m_{\bpi_1}.
$$
In particular, if we denote by $M_0$ the $\bu(L^e(\lie h))$-submodule
of $M$ generated by $m_{\bpi_1}$, then we have,
$$ M_0=\bigoplus_{r\in\bz} M_{\lambda_{\bpi_1}+r\delta}\cong \bigoplus_{r\in\bz} \left((V_1)_{\lambda_{\bpi_1}+r\delta}\oplus
(V_2)_{\lambda_{\bpi_1}+r\delta}\right)
$$
as $\lie h_e$-modules.

If  $\dim M_{\lambda_{\bpi_1}+r\delta}\le 1$ (in particular, if  $\lambda_{\bpi_1}\ne \lambda_{\bpi_2}$), then
$M$ is an $\ell$-highest weight module for
$L^e(\lie g)$. Then by~\propref{CI80}  $M$ is a quotient of
$L^s(\mathcal W(\bpi_1))$. Since ~\propref{L60} implies that $L^s(\mathcal W(\bpi_1))\in \operatorname{Ob}\mathcal I^{fin}_{\bar\chi_{\bpi_1}}$ it follows that $M\in \operatorname{Ob}\mathcal I^{fin}_{\bar\chi_{\bpi_1}}$. That is impossible if
 $\bar\chi_{\bpi_2}\ne\bar\chi_{\bpi_1}$ and we obtain a contradiction.

If $ \lambda_{\bpi_1}=\lambda_{\bpi_2}=\lambda$ and $\dim M_{\lambda+p_0\delta}=2$ for some $p_0\in\bz$, then
$\dim M_{\lambda+p\delta}=2$ for all $p=p_0\pmod k$ where $k={\operatorname{lcm}(m(\bpi_1),m(\bpi_2))}$.
Fix~$m_{\bpi_1,l}\in M_{\lambda+(p_0+kl)\delta}$ such that~$\xi(m_{\bpi_1,l})=v_{\bpi_1}\tensor t^{p_0+kl}$ and set~$m_{\bpi_2,l}=
\phi(v_{\bpi_2}\tensor t^{p_0+kl}) \in
M_{\lambda+(p_0+kl)\delta}$. Then~$\{m_{\bpi_1,l}$, $m_{\bpi_2,l}\}$ is a basis of~$M_{\lambda+(p_0+kl)\delta}$.
Let $\bar\psi_\bpi$ denote the composition of~$\psi_\bpi$ with the map~$\bc[t^{\pm1}]\to\bc$, $t\mapsto1$.
Since~$\bar\chi_{\bpi_1}\not=\bar\chi_{\bpi_2}$ there exists $l,s\in\bz$ and~$x\in\bu(L(\lie h))_{ks}$
such that $\bar\psi_{\bpi_1}(x)\not=\bar\psi_{\bpi_2}(x)$ and
$$x m_{\bpi_1,l}=\bar\psi_{\bpi_1}(x) m_{\bpi_1,l+s}+c m_{\bpi_2,l+s},\qquad x m_{\bpi_2,l}=\bar\psi_{\bpi_2}(x) m_{\bpi_2,l+s},$$
where~$c\in\bc$.
Let~$\gamma: M_{\lambda+(p_0+(l+s)k)\delta}\to  M_{\lambda+(p_0+kl)\delta}$ be the isomorphism of
vector spaces defined by~$\gamma(m_{\bpi_j,l+s})=m_{\bpi_j,l}$, $j=1,2$. Then the matrix of~$\gamma\circ x$ in the basis~$\{m_{\bpi_1,l}, m_{\bpi_2,l}\}$ is upper triangular
with distinct diagonal entries and so $\gamma\circ x$ has two one-dimensional eigenspaces.
Obviously, $m_{\bpi_2,l}$ is an eigenvector of $\gamma\circ x$ corresponding to the eigenvalue $\bar\psi_{\bpi_2}(x)$.
On the other hand, since~$m_{\bpi_2,l}\in \bu(L^e(\lie h))m_{\bpi_1,l}$, it follows that~$m_{\bpi_2,l}=y m_{\bpi_1,l}$
for some~$y\in\bu(L^e(\lie h))_0$. Now
$$
y\gamma(x m_{\bpi_1,l})=\gamma(x y m_{\bpi_1,l}),$$ and hence $m_{\bpi_2,l}$ is also an eigenvector of~$\gamma\circ x$ corresponding to
the eigenvalue~$\bar\psi_{\bpi_1}(x)$, which is a contradiction.
\end{pf}

\begin{cor}  Suppose that~$V_i\in\operatorname{Ob}\mathcal I^{fin}_{\bar\chi_i}$, $i=1,2$ with~$\bar\chi_1\not=\bar\chi_2$.
Then $\Ext^1_{\mathcal I^{fin}}(V_1,V_2)=0$.
\end{cor}
\begin{pf} An  induction on the pseudo-length (similar to the induction on length used in
\cite[Lemma~5.2]{CM} for $\cal F$) yields the corollary.
\end{pf}

\subsection{}\label{CI120}
We are now  able to prove~\propref{CI20}.
\begin{pf} Let~$V$ is an indecomposable object in~$\mathcal
I^{fin}$. We proceed by induction on the pseudo-length of~$V$.
If~$V$ is simple then it clear from the definition
that~$V\in\mathcal I^{fin}_{\bar\chi_\bpi}$ for
some~$\bpi\in\mathcal P^+$. Otherwise, we have an extension
\begin{equation}\label{CI120.10}
\begin{CD}
0 @>>> V_1 @>>> V @>>> U @>>> 0,
\end{CD}
\end{equation}
where either~$V_1=V^{L(\lie g)}$ or~$V_1\not\subset V^{L(\lie g)}$ and is simple. In any case,
$V_1$ is an object in~$\mathcal I^{fin}_{\bar\chi}$ for some
$\chi\in\Xi$. Using~\thmref{main1}, we can
write~$$U=\bigoplus_{j=0}^k U_j$$ where~$U_0=\bigoplus_{r\in\bz}
\bc_{r\delta}^{m_r}$, $m_r\ge 0$, and for $j=1,\dots ,k$ the
module $U_j$ is indecomposable with ${\wt}^e(U_j)\not\subset
\bz\delta$. Notice that the pseudo-length of~$U_j$, $j\ge 1$ is
strictly smaller than that of~$V$ unless $U$ is indecomposable and
$V_1=V^{L(\lie g)}$.

Suppose first that $V_1\ne V^{L(\lie g)}$. Then,  by the induction
hypothesis, we see that for $j=1,\dots ,k$,  $U_j$ is an object
in~$\mathcal I^{fin}_{\bar\chi_j}$ for some~$\chi_j\in\Xi$.
Suppose that $\bar\chi\not=\bar\chi_{j_0}$ for some $0\le j_0\le
k$. Then by~\corref{CI100}
$$
\Ext^1_{\mathcal I^{fin}}(U,V_1)\cong\bigoplus_{j\not=j_0}
\Ext^1_{\mathcal I^{fin}}(U_j,V_1).
$$
In other words, the exact sequence~\eqref{CI120.10} is equivalent to
$$
\begin{CD}
0 @>>> V_1 @>>> U_{j_0} \oplus V' @>>> U_{j_0}\oplus
\bigoplus_{j\not=j_0} U_j @>>> 0,
\end{CD}
$$
where
$$
\begin{CD}
0 @>>> V_1 @>>> V' @>>> \bigoplus_{j\not=j_0} U_j @>>> 0
\end{CD}
$$
is an element of $\bigoplus_{j\not=j_0} \Ext^1_{\mathcal
I^{fin}}(U_j,V_1)$. This is a contradiction since~$V$ is
indecomposable. Finally, suppose that~$U$ is indecomposable
and~$V_1=V^{L(\lie g)}$. Then~$U^{L(\lie g)}=0$ and so we have a
short exact sequence
$$
\begin{CD}
0 @>>> U_1 @>>> U @>>> U_2 @>>>0,
\end{CD}
$$
where~$U_1$ is simple and $U_1\not\subset U^{L(\lie g)}$. Hence the
pseudo-length of~$U_2$ is strictly smaller than that of~$U$, and so we
can apply the above argument and show that~$U$ is an object
in~$\mathcal I^{fin}_{\bar\chi}$ for some~$\chi\in\Xi$.
If~$\bar\chi\not=0$, then $\Ext^1_{\mathcal I^{fin}}(U,V_1)=0$ by~\corref{CI100} and
so~$V\cong V_1\oplus U$ which is a contradiction. Thus, $\bar\chi=0$ and~$V$ is an object in~$\operatorname{Ob} \mathcal I^{fin}_{\bar 0}$.
\end{pf}

\section*{List of notations}
\def\bqq{{\setbox0\hbox{$\wh{\bu}_q^+$}\setbox2\null\ht2\ht0\dp2\dp0\box2}}
\noindent
\begin{tabular}{p{1.8in}@{\bqq}l}
$\lie g$, $\lie h$ &\ref{P10}, p.~\pageref{P10}\\
$I$, $\alpha_i$, $\varpi_i$&\ref{P10}, p.~\pageref{P10}\\
$R$, $R^+$, $Q$, $Q^+$, $P$, $P^+$&\ref{P10}, p.~\pageref{P10}\\
$\theta$, $\theta_s$&\ref{P10}, p.~\pageref{P10}\\
$W$&\ref{P10}, p.~\pageref{P10}\\
$x_{\alpha}^\pm$, $\alpha^\vee$&\ref{P10}, p.~\pageref{P10}\\
$\Gamma$, $\varpi_\gamma$&\ref{P10}, p.~\pageref{P10}\\
$V(\lambda)$, $v_\lambda$&\ref{P15}, p.~\pageref{P15}\\
$\bu(\lie a)$&\ref{P20}, p.~\pageref{P20}\\
$L(\lie a)$, $L^e(\lie a)$, $d$&\ref{P20}, p.~\pageref{P20}\\
$\lie h_e$&\ref{P25}, p.~\pageref{P25}\\
$\delta$&\ref{P25}, p.~\pageref{P25}\\
$P_e$, $P_e^+$&\ref{P25}, p.~\pageref{P25}\\
$\wh{W}$&\ref{P25}, p.~\pageref{P25}\\
$r_\lambda$&\ref{P25}, p.~\pageref{P25}\\
$\Lambda_i^\pm(u)$, $\Lambda_{i,\pm k}$&\ref{P30}, p.~\pageref{P30}\\
$\cal I$&\ref{E20}, p.~\pageref{E20}\\
$\wt^e(V)$&\ref{E20}, p.~\pageref{E20}\\
$\bc_{a\delta}$&\ref{E40}, p.~\pageref{E40}
\end{tabular}
\hfill
\begin{tabular}{p{1.8in}@{\bqq}l}
$\cal I^{fin}$&\ref{E40}, p.~\pageref{E40}\\
$V^+_\lambda$&\ref{E60}, p.~\pageref{E60}\\
$V^{L(\lie g)}$&\ref{E100}, p.~\pageref{E100}\\
$(\cdot,\cdot)_{\lie g}$&\ref{Fin}, p.~\pageref{Fin}\\
$\sim$&\ref{F120}, p.~\pageref{F120}\\
$\Xi$, $\overline{\Xi}$&\ref{F140}, p.~\pageref{F140}\\
$\chi_{i,a}$&\ref{F140}, p.~\pageref{F140}\\
$\cal P^+$&\ref{L10}, p.~\pageref{L10}\\
$\bvpi_{i,a}$, $\bpi_{\lambda,a}$&\ref{L10}, p.~\pageref{L10}\\
$\lambda_\bpi$, $m(\bpi)$&\ref{L10}, p.~\pageref{L10}\\
$\cal F$&\ref{L20}, p.~\pageref{L20}\\
$\mathcal W(\bpi)$, $w_\bpi$, $\mathcal V(\bpi)$, $v_\bpi$&\ref{L20}, p.~\pageref{L20}\\
$\tau_a$&\ref{L20}, p.~\pageref{L20}\\
$\chi_\bpi$&\ref{L40}, p.~\pageref{L40}\\
$\cal F_\chi$&\ref{L40}, p.~\pageref{L40}\\
$L$&\ref{L60}, p.~\pageref{L60}\\
$L^s(\mathcal V(\bpi))$&\ref{L60}, p.~\pageref{L60}\\
$\cal I^{fin}_{\overline{\chi}}$&\ref{L80}, p.~\pageref{L80}\\
${}^\#$&\ref{L90}, p.~\pageref{L90}\\
\end{tabular}

\bibliographystyle{amsplain}

\end{document}